\documentclass[12pt]{article}
\usepackage{amssymb}
\usepackage{amsfonts}

\begin{document}

\title{Pricing of basket options II}
\author{Alexander Kushpel \\
Department of Mathematics,\\
University of Leicester, LE1 7RH\\
E-mail: ak412@le.ac.uk}
\date{31 January 2014}
\maketitle

\begin{abstract}
We consider the problem of approximation of density functions which is
important in the theory of pricing of basket options. Our method is well
adopted to the multidimensional case. Observe that implementations of polynomial
and spline approximation in this situation are connected with difficulties
of fundamental nature. A simple approximation formula for European call
options is presented. It is shown that this approximation formula has
exponential rate of convergence.
\end{abstract}



\section{Approximation of density functions}

\bigskip

We shall use notations adopted in \cite{kushpel-basket2}. Remind that the
pricing formula has the form $V=e^{-rT}\mathbb{E}^{\mathbb{Q}}\left[ \varphi %
\right] $. Since the reward function $\varphi $ has a simple structure the
main problem is to approximate the respective risk-neutral density function.
For our model the characteristic function $\Phi ^{\mathbb{Q}}\left( \mathbf{%
v,}T\right) $ has the form
\begin{equation}
\Phi ^{\mathbb{Q}}\left( \mathbf{v,}T\right) =\prod\limits_{s=1}^{n}\exp
\left( -T\psi _{s}^{\mathbb{Q}}\left( v_{s}\right) \right) \cdot
\prod\limits_{m=1}^{n}\exp \left( -T\phi _{m}^{\mathbb{Q}}\left(
\sum_{k=1}^{n}a_{k,m}v_{k}\right) \right) ,  \label{phiqvt}
\end{equation}%
and the respective characteristic exponent
\[
\tau ^{\mathbb{Q}}\left( \mathbf{v,}T\right) :=-T\left( \sum_{s=1}^{n}\psi
_{s}^{\mathbb{Q}}\left( v_{s}\right) +\sum_{m=1}^{n}\phi _{m}^{\mathbb{Q}%
}\left( \sum_{k=1}^{n}a_{k,m}v_{k}\right) \right)
\]%
satisfies the condition%
\begin{equation}
\tau ^{\mathbb{Q}}\left( -i\mathbf{e}_{s},T\right) =-r,1\leq s\leq n.
\label{EMM condition}
\end{equation}%
In particular, in the case of KoBoL family
\[
\psi _{s}^{\mathbb{Q}}\left( \xi _{s}\right) =-i\mu _{s}\xi
_{s}+c_{s,+}\Gamma \left( -\nu _{s}\right) \left( \left( -\lambda
_{s,-}\right) ^{\nu _{s}}-\left( -\lambda _{s,-}-i\xi _{s}\right) ^{\nu
_{s}}\right)
\]%
\[
+c_{s,-}\Gamma \left( -\nu _{s}\right) \left( \lambda _{s,+}^{\nu
_{s}}-\left( \lambda _{s,+}+i\xi _{s}\right) ^{\nu _{s}}\right) ,
\]%
\begin{equation}
\lambda _{s,-}<0<\lambda _{s,+},\nu _{s}\in \left( 0,1\right) ,\mu _{s}\in
\mathbb{R},c_{s,+}>0,c_{s,-}>0,1\leq s\leq n.  \label{s-exponent}
\end{equation}%
It is possible to show \cite{bl1} that for any $\nu \in \left( 0,1\right)
\cup \left( 1,2\right) ,$%
\[
\psi ^{\mathbb{Q}}\left( \xi \right) =-i\mu \xi +c_{+}\Gamma \left( -\nu
\right) \left( \left( -\lambda _{-}\right) ^{\nu }-\left( -\lambda _{-}-i\xi
\right) ^{\nu }\right)
\]%
\[
+c_{-}\Gamma \left( -\nu \right) \left( \lambda _{+}^{\nu }-\left( \lambda
_{+}+i\xi \right) ^{\nu }\right)
\]%
and
\[
\psi ^{\mathbb{Q}}\left( \xi \right) =-i\mu \xi +c_{+}\left( \ln \left(
-\lambda _{-}-i\xi \right) -\ln \left( -\lambda _{-}\right) \right)
\]%
\[
+c_{-}\left( \ln \left( \lambda _{+}-i\xi \right) -\ln \lambda _{+}\right)
,\nu =0,
\]%
\[
\psi ^{\mathbb{Q}}\left( \xi \right) =-i\mu \xi +c_{+}\left( \left( -\lambda
_{-}\right) \ln \left( -\lambda _{-}\right) -\left( -\lambda _{-}-i\xi
\right) \ln \left( -\lambda _{-}-i\xi \right) \right)
\]%
\[
+c_{-}\left( \lambda _{+}\ln \lambda _{+}-\left( \lambda _{+}-i\xi \right)
\ln \left( \lambda _{+}-i\xi \right) \right) ,\nu =1.
\]
Assume that all characteristic exponents $\psi _{s}^{\mathbb{Q}},1\leq s\leq
n$ and $\phi _{m}^{\mathbb{Q}},1\leq m\leq n$ correspond a KoBoL process and
hence analytically extendable into the strips $\lambda _{s,-}<\kappa
_{s,-}<0<\kappa _{s,+}<\lambda _{s,+}$, $1\leq s\leq n$ and $\lambda
_{s,-}^{^{\prime }}<\kappa _{s,-}^{^{\prime }}<0<\kappa _{s,+}^{^{\prime
}}<\lambda _{s,+}^{^{\prime }}$, $1\leq s\leq n$ respectively. In this case $%
\Phi ^{\mathbb{Q}}\left( \mathbf{v,}t\right) =\Phi ^{\mathbb{Q}}\left(
v_{1},\cdot \cdot \cdot ,v_{n}\mathbf{,}t\right) $ admits an analytic
extension into the tube $T_{n}\subset \mathbb{C}^{n},$
\[
T_{n}:=\left( \prod\limits_{s=1}^{n}\left\{ \Im v_{s}\in \left[ \kappa
_{s,-},\kappa _{s,+}\right] \right\} \right)
\]%
\[
\cap \left( \left\{ \Im \left( \sum_{k=1}^{n}a_{k,m}v_{k}\right) \in \left[
\kappa _{m,-}^{^{\prime }},\kappa _{m,+}^{^{\prime }}\right] ,1\leq m\leq
n\right\} \right) .
\]%
Let%
\[
a_{+}:=\min_{1\leq s,m\leq n}\left\{ \kappa _{s,+},\kappa _{m,+}^{^{\prime
}}\min \left\{ 1,\left( \sum_{k=1}^{n}a_{k,m}\right) ^{-1}\right\} \right\}
\]%
and%
\[
a_{-}:=\max_{1\leq s,m\leq n}\left\{ \kappa _{s,-},\kappa _{m,-}^{^{\prime
}}\min \left\{ 1,\left( \sum_{k=1}^{n}a_{k,m}\right) ^{-1}\right\} \right\}
\]%
then
\[
T_{n}^{^{\prime }}:=\prod\limits_{s=1}^{n}\left\{ \Im v_{s}\in \lbrack
a_{-},a_{+}]\right\} \subset T_{n}.
\]
Let $\mathbf{a}_{+}:=\left( a_{+},\cdot \cdot \cdot ,a_{+}\right) $ and $%
\mathbf{a}_{-}:=\left( a_{-},\cdot \cdot \cdot ,a_{-}\right) $.

{\bf Theorem 1.} \label{representation2014} 
\begin{em}
Let $\psi _{s}^{\mathbb{Q}},1\leq s\leq n$ and $%
\phi _{m}^{\mathbb{Q}},1\leq m\leq n$ defined by (\ref{s-exponent}) then the
respective density function $p_{t}^{\mathbb{Q}}\left( \cdot \right) $ can be
represented as
\[
p_{T}^{\mathbb{Q}}\left( \cdot \right) =\frac{\left( 2\pi \right) ^{-n}}{%
\exp \left( \left\langle \cdot ,\mathbf{a}_{+}\right\rangle \right) +\exp
\left( \left\langle \cdot ,\mathbf{a}_{-}\right\rangle \right) }
\]%
\[
\times \int_{\mathbb{R}^{n}}\exp \left( -i\left\langle \cdot \mathbf{,v}%
\right\rangle \right) \left( \Phi ^{\mathbb{Q}}\left( \mathbf{v-}i\mathbf{a}%
_{+}\mathbf{,}T\right) +\Phi ^{\mathbb{Q}}\left( \mathbf{v-}i\mathbf{a}_{-}%
\mathbf{,}T\right) \right) d\mathbf{v}.
\]%
In particular, if $-\mathbf{a}_{-}=\mathbf{a}_{+}:=\mathbf{a}$ then
\[
p_{T}^{\mathbb{Q}}\left( \cdot \right) =\frac{1}{2\left( 2\pi \right) ^{n}}%
\left( \cosh \left( \left\langle \cdot ,\mathbf{a}\right\rangle \right)
\right) ^{-1}
\]%
\[
\times \int_{\mathbb{R}^{n}}\exp \left( -i\left\langle \cdot \mathbf{,v}%
\right\rangle \right) \left( \Phi ^{\mathbb{Q}}\left( \mathbf{v+}i\mathbf{a,}%
T\right) +\Phi ^{\mathbb{Q}}\left( \mathbf{v-}i\mathbf{a,}T\right) \right) d%
\mathbf{v}.
\]
\end{em}

{\it Proof}
In our notations the density function can be represented as
\[
p_{T}^{\mathbb{Q}}\left( \cdot \right) =\left( 2\pi \right) ^{-n}\int_{%
\mathbb{R}^{n}}\exp \left( -i\left\langle \cdot \mathbf{,v}\right\rangle
\right) \Phi ^{\mathbb{Q}}\left( \mathbf{v,}T\right) d\mathbf{v}=\left( 2\pi
\right) ^{-n}\mathbf{F}\left( \Phi ^{\mathbb{Q}}\left( \mathbf{v,}T\right)
\right) \left( \cdot \right) .
\]%
It is easy to check that $\psi _{s}^{\mathbb{Q}}\left( \xi _{s}\right) $, $%
1\leq s\leq n$ admits an analytic extension into the strip $\Im \xi _{s}\in %
\left[ \kappa _{s,-},\kappa _{s,+}\right] $, where $\lambda _{s,-}<\kappa
_{s,-}<0<\kappa _{s,+}<\lambda _{s,+}$, $1\leq s\leq n$ and
\[
\lim_{\Im \xi _{s}\in \left[ \kappa _{s,-},\kappa _{s,+}\right] ,\left\vert
\xi _{s}\right\vert \rightarrow \infty }\frac{\psi _{s}^{\mathbb{Q}}\left(
\xi _{s}\right) }{\eta _{s}\left( \left\vert \Re \xi _{s}\right\vert \right)
}=1,
\]%
\begin{equation}
\eta _{s}\left( x\right) =\left\{
\begin{array}{cc}
-i\mu _{s}x-\Gamma \left( -\nu _{s}\right) \left( c_{s,+}\exp \left( \frac{%
-i\nu _{s}\pi }{2}\right) +c_{s,-}\exp \left( \frac{i\nu _{s}\pi }{2}\right)
\right) x^{\nu }, & \nu _{s}\in \left( 0,1\right) \cup \left( 1,2\right) ,
\\
-i\left( \mu _{s}x+c_{s,+}x\ln x-c_{s,-}x\ln x\right) +\frac{\pi \left(
c_{s,+}+c_{s,-}\right) x}{2}, & \nu _{s}=1, \\
-i\mu _{s}x+\left( c_{s,+}+c_{s,-}\right) \ln x, & \nu _{s}=0.%
\end{array}%
\right.  \label{asymp 1}
\end{equation}%
In particular, if $c_{s,+}=c_{s,-}=c_{s}$ then $\eta _{s}\left( x\right) $
simplifies as
\begin{equation}
\eta _{s}\left( x\right) =\left\{
\begin{array}{cc}
-i\mu _{s}x-2c_{s}\Gamma \left( -\nu _{s}\right) \cos \left( \frac{\nu
_{s}\pi }{2}\right) x^{\nu }, & \nu _{s}\in \left( 0,1\right) \cup \left(
1,2\right) , \\
-i\mu _{s}x+\pi c_{s}x, & \nu =1, \\
-i\mu _{s}x+2\ln x, & \nu =0.%
\end{array}%
\right.  \label{asymp 2}
\end{equation}%
The same asymptotic are valid for $\phi _{m}^{\mathbb{Q}},1\leq m\leq n.$
Observe that
\[
\mathrm{sign}\left( -\Gamma \left( -\nu \right) \cos \frac{\pi \nu }{2}%
\right) >0.
\]%
Then from (\ref{asymp 1}) and (\ref{asymp 2}) it follows that
\[
\lim_{\left\langle \mathbf{v},\mathbf{v}\right\rangle \rightarrow \infty ,%
\mathbf{v}\in T_{n}^{^{\prime }}}\Phi ^{\mathbb{Q}}\left( \mathbf{v,}%
T\right) d\mathbf{v}=0.
\]%
Hence, applying Cauchy's theorem $n$ times in the tube $T_{n}^{^{\prime }}$,
which is justified by (\ref{asymp 1}), we get
\[
p_{T}^{\mathbb{Q}}\left( \cdot \right) =\left( 2\pi \right) ^{-n}\int_{%
\mathbb{R}^{n}}\exp \left( -i\left\langle \cdot \mathbf{,v}\right\rangle
\right) \Phi ^{\mathbb{Q}}\left( \mathbf{v,}T\right) d\mathbf{v}
\]%
\[
=\left( 2\pi \right) ^{-n}\int_{\mathbb{R}^{n}+i\mathbf{a}_{+}}\exp \left(
-i\left\langle \cdot \mathbf{,v}\right\rangle \right) \Phi ^{\mathbb{Q}%
}\left( \mathbf{v,}T\right) d\mathbf{v}
\]%
\[
=\left( 2\pi \right) ^{-n}\int_{\mathbb{R}^{n}}\exp \left( -i\left\langle
\cdot \mathbf{,x-}i\mathbf{a}_{+}\right\rangle \right) \Phi ^{\mathbb{Q}%
}\left( \mathbf{x-}i\mathbf{a}_{+}\mathbf{,}T\right) d\mathbf{x}
\]%
\[
=\exp \left( -\left\langle \cdot ,\mathbf{a}_{+}\right\rangle \right) \left(
2\pi \right) ^{-n}\int_{\mathbb{R}^{n}}\exp \left( -i\left\langle \cdot
\mathbf{,x}\right\rangle \right) \Phi ^{\mathbb{Q}}\left( \mathbf{x-}i%
\mathbf{a}_{+}\mathbf{,}T\right) d\mathbf{x,}
\]%
or%
\begin{equation}
p_{T}^{\mathbb{Q}}\left( \cdot \right) \exp \left( \left\langle \cdot ,%
\mathbf{a}_{+}\right\rangle \right) =\left( 2\pi \right) ^{-n}\int_{\mathbb{R%
}^{n}}\exp \left( -i\left\langle \cdot \mathbf{,x}\right\rangle \right) \Phi
^{\mathbb{Q}}\left( \mathbf{x-}i\mathbf{a}_{+}\mathbf{,}T\right) d\mathbf{x.}
\label{ptq1}
\end{equation}%
Similarly,
\begin{equation}
p_{T}^{\mathbb{Q}}\left( \cdot \right) \exp \left( \left\langle \cdot ,%
\mathbf{a}_{-}\right\rangle \right) =\left( 2\pi \right) ^{-n}\int_{\mathbb{R%
}^{n}}\exp \left( -i\left\langle \cdot \mathbf{,x}\right\rangle \right) \Phi
^{\mathbb{Q}}\left( \mathbf{x-}i\mathbf{a}_{-}\mathbf{,}T\right) d\mathbf{x.}
\label{ptq2}
\end{equation}%
Comparing (\ref{ptq1}) and (\ref{ptq2}) we get the proof.

We will need the following result which is known as the Poisson summation
formula.

{\bf Theorem 2.}
\label{Poisson summation} 
\begin{em} (\cite{stain} p.252) Suppose that for some $A>0$
and $\delta >0$ we have
\[
\max \left\{ f\left( \mathbf{x}\right) ,\mathbf{F}f\left( \mathbf{x}\right)
\right\} \leq A\left( 1+\left\vert \mathbf{x}\right\vert \right) ^{-n-\delta
}
\]%
then
\[
\sum_{\mathbf{m}\in \mathbb{Z}^{n}}f\left( \mathbf{x}+P\mathbf{m}\right) =%
\frac{1}{P^{n}}\sum_{\mathbf{m}\in \mathbb{Z}^{n}}\mathbf{F}f\left( \frac{%
\mathbf{m}}{P}\right) \exp \left( i\left\langle \frac{2\pi \mathbf{m}}{P},%
\mathbf{x}\right\rangle \right)
\]%
for any $P>0$. The series converges converges absolutely.
\end{em}

Put
\[
M_{T}=M_{T}\left( \Phi ^{\mathbb{Q}},\mathbf{a}\right)
\]%
\[
:=\frac{1}{2\left( 2\pi \right) ^{n}}\left\Vert \int_{\mathbb{R}^{n}}\exp
\left( -i\left\langle \cdot \mathbf{,v}\right\rangle \right) \left( \Phi ^{%
\mathbb{Q}}\left( \mathbf{v+}i\mathbf{a,}T\right) +\Phi ^{\mathbb{Q}}\left(
\mathbf{v-}i\mathbf{a,}T\right) \right) d\mathbf{v}\right\Vert _{\infty },
\]%
fix $T>0,\epsilon >0$ and select such $P\in \mathbb{N}$ that
\begin{equation}
M_{T}\sum_{\mathbf{m}\in \mathbb{Z}^{n}\diagdown \left\{ \mathbf{0}\right\}
}\left( \cosh \left( \left\langle \frac{2P-1}{2}\mathbf{m},\mathbf{a}%
\right\rangle \right) \right) ^{-1}\leq \epsilon .  \label{epsilon111}
\end{equation}

{\bf Theorem 3.} 
\label{desity 111} 
\begin{em}
Let $-\mathbf{a}_{-}=\mathbf{a}_{+}:=\mathbf{a}$ then in
our notations%
\[
\left\Vert p_{T}^{\mathbb{Q}}\left( \mathbf{x}\right) -\frac{1}{P^{n}}\sum_{%
\mathbf{m}\in \mathbb{Z}^{n}}\Phi ^{\mathbb{Q}}\left( -\frac{\mathbf{m}}{P}%
,T\right) \exp \left( i\left\langle \frac{2\pi \mathbf{m}}{P},\mathbf{x}%
\right\rangle \right) \right\Vert _{L_{p}\left( \frac{P}{2}Q_{n}\right)
}\leq \epsilon P^{n/p},
\]%
where $1\leq p\leq \infty .$
\end{em}

{\bf Proof}
Using Theorem \ref{representation2014} we get
\begin{equation}
p_{T}^{\mathbb{Q}}\left( \cdot \right) \leq M_{T}\left( \cosh \left(
\left\langle \cdot ,\mathbf{a}\right\rangle \right) \right) ^{-1}.
\label{majorant111}
\end{equation}%
Applying (\ref{majorant111}) we can check that the conditions of Theorem \ref%
{Poisson summation} are satisfied. Hence using condition (\ref{epsilon111})
we get%
\[
\epsilon \geq \left\Vert \sum_{\mathbf{m}\in \mathbb{Z}^{n}\diagdown \left\{
\mathbf{0}\right\} }p_{T}^{\mathbb{Q}}\left( \mathbf{x}+P\mathbf{m}\right)
\right\Vert _{L_{\infty }\left( \frac{P}{2}Q_{n}\right) }
\]%
\[
=\left\Vert p_{T}^{\mathbb{Q}}\left( \mathbf{x}\right) -\sum_{\mathbf{m}\in
\mathbb{Z}^{n}}p_{T}^{\mathbb{Q}}\left( \mathbf{x}+P\mathbf{m}\right)
\right\Vert _{L_{\infty }\left( \frac{P}{2}Q_{n}\right) }
\]%
\[
=\left\Vert p_{T}^{\mathbb{Q}}\left( \mathbf{x}\right) -\frac{1}{P^{n}}\sum_{%
\mathbf{m}\in \mathbb{Z}^{n}}\mathbf{F}p_{T}^{\mathbb{Q}}\left( \frac{%
\mathbf{m}}{P}\right) \exp \left( i\left\langle \frac{2\pi \mathbf{m}}{P},%
\mathbf{x}\right\rangle \right) \right\Vert _{L_{\infty }\left( \frac{P}{2}%
Q_{n}\right) },
\]%
where $Q_{n}:=\left\{ \mathbf{x}\left\vert \mathbf{x}=\left( x_{1},\cdot
\cdot \cdot ,x_{n}\right) \in \mathbb{R}^{n},\left\vert x_{k}\right\vert
\leq 1,1\leq k\leq n\right. \right\} $ is the unit qube in $\mathbb{R}^{n}$.
Observe that
\[
\Phi ^{\mathbb{Q}}\left( -\mathbf{x,}T\right) =\left( 2\pi \right) ^{n}%
\mathbf{F}^{-1}p_{T}^{\mathbb{Q}}\left( -\mathbf{x}\right)
\]%
\[
=\left( 2\pi \right) ^{n}\left( \frac{1}{\left( 2\pi \right) ^{n}}\right)
\int_{\mathbb{R}^{n}}\mathbf{\exp }\left( i\left\langle \mathbf{x,y}%
\right\rangle \right) p_{T}^{\mathbb{Q}}\left( -\mathbf{y}\right) d\mathbf{y}
\]%
\[
=\int_{\mathbb{R}^{n}}\mathbf{\exp }\left( i\left\langle -\mathbf{x,y}%
\right\rangle \right) p_{T}^{\mathbb{Q}}\left( \mathbf{y}\right) d\mathbf{y}
\]%
\[
=\mathbf{F}p_{T}^{\mathbb{Q}}\left( \mathbf{x}\right) .
\]%
Consequently,
\[
\left\Vert p_{T}^{\mathbb{Q}}\left( \mathbf{x}\right) -\frac{1}{P^{n}}\sum_{%
\mathbf{m}\in \mathbb{Z}^{n}}\Phi ^{\mathbb{Q}}\left( -\frac{\mathbf{m}}{P}%
,T\right) \exp \left( i\left\langle \frac{2\pi \mathbf{m}}{P},\mathbf{x}%
\right\rangle \right) \right\Vert _{L_{\infty }\left( \frac{P}{2}%
Q_{n}\right) }\leq \epsilon
\]%
and%
\[
\left\Vert p_{T}^{\mathbb{Q}}\left( \mathbf{x}\right) -\frac{1}{P^{n}}\sum_{%
\mathbf{m}\in \mathbb{Z}^{n}}\Phi ^{\mathbb{Q}}\left( -\frac{\mathbf{m}}{P}%
,T\right) \exp \left( i\left\langle \frac{2\pi \mathbf{m}}{P},\mathbf{x}%
\right\rangle \right) \right\Vert _{L_{1}\left( \frac{P}{2}Q_{n}\right)
}\leq \epsilon P^{n}.
\]%
Finally, applying Riesz-Thorin interpolation theorem we obtain
\[
\left\Vert p_{T}^{\mathbb{Q}}\left( \mathbf{x}\right) -\frac{1}{P^{n}}\sum_{%
\mathbf{m}\in \mathbb{Z}^{n}}\Phi ^{\mathbb{Q}}\left( -\frac{\mathbf{m}}{P}%
,T\right) \exp \left( i\left\langle \frac{2\pi \mathbf{m}}{P},\mathbf{x}%
\right\rangle \right) \right\Vert _{L_{p}\left( \frac{P}{2}Q_{n}\right)
}\leq \epsilon P^{n/p}.
\]

Observe that the function $\left\vert \Phi ^{\mathbb{Q}}\left( -\frac{%
\mathbf{m}}{P},T\right) \right\vert $ exponentially decays as $\left\vert
\mathbf{m}\right\vert \rightarrow \infty $. Hence the series
\[
\frac{1}{P^{n}}\sum_{\mathbf{m}\in \mathbb{Z}^{n}}\Phi ^{\mathbb{Q}}\left( -%
\frac{\mathbf{m}}{P},T\right) \exp \left( i\left\langle \frac{2\pi \mathbf{m}%
}{P},\mathbf{x}\right\rangle \right)
\]%
converges absolutely and represents an infinitely differentiable function.
The next statements deal with approximation of this function using $M$-term
exponential sums.

{\bf Theorem 4.}
\label{m-term1} 
\begin{em}
Let
\[
\Omega \left( \Phi ^{\mathbb{Q}}\left( \mathbf{v,}T\right) ,\varrho \right)
:=\left\{ \mathbf{v}\in \mathbb{R}^{n},\left\vert \Phi ^{\mathbb{Q}}\left(
\mathbf{v,}T\right) \right\vert \geq \varrho \right\} ,
\]%
\[
\varkappa _{n}=2^{n}\Gamma ^{-1}\left( 1+\sum_{k=1}^{n}\nu _{k}\right)
\left( \prod\limits_{s=1}^{n}\Gamma \left( 1+\nu _{s}\right) \right) \left(
\prod\limits_{s=1}^{n}\left( d_{s}T\right) ^{-\nu _{s}^{-1}}\right)
\]%
and
\[
d_{s}:=-\Gamma \left( -\nu _{s}\right) \cos \left( \frac{\nu _{s}\pi }{2}%
\right) \left( c_{s,+}+c_{s,-}\right) ,1\leq s\leq n.
\]%
Then
\[
\mathrm{Card}\left( \Omega \left( \Phi ^{\mathbb{Q}}\left( \mathbf{v,}%
T\right) ,\varrho \right) \cap \mathbb{Z}^{n}\right) \lesssim \varkappa
_{n}\left( \ln \varrho ^{-1}\right) ^{\sum_{k=1}^{n}\nu _{s}^{-1}},\varrho
\rightarrow 0.
\]
\end{em}

{\bf Proof}
Observe that
\[
p_{T}^{\mathbb{Q}}\left( \cdot \right) =\left( 2\pi \right) ^{-n}\mathbf{F}%
\left( \Phi ^{\mathbb{Q}}\left( \mathbf{v,}T\right) \right) \left( \cdot
\right) .
\]%
From (\ref{asymp 1}) it follows that $\Phi ^{\mathbb{Q}}\in C\left( \mathbb{R%
}^{n}\right) \cap L_{2}\left( \mathbb{R}^{n}\right) $. Hence by Plancherel's
theorem
\[
\mathbf{F}^{-1}\left( p_{T}^{\mathbb{Q}}\right) \left( \cdot \right) =\left(
2\pi \right) ^{-n}\Phi ^{\mathbb{Q}}\left( \mathbf{v,}T\right) ,
\]%
or%
\[
\Phi ^{\mathbb{Q}}\left( \mathbf{v,}T\right) =\int_{\mathbb{R}^{n}}\exp
\left( i\left\langle \mathbf{v,x}\right\rangle \right) p_{T}^{\mathbb{Q}%
}\left( \mathbf{x}\right) d\mathbf{x}.
\]%
Since
\[
p_{T}^{\mathbb{Q}}\left( \mathbf{x}\right) \geq 0,\int_{\mathbb{R}%
^{n}}p_{T}^{\mathbb{Q}}\left( \mathbf{x}\right) d\mathbf{x=1}
\]%
then
\begin{equation}
\max \left\{ \mathbf{v}\in \mathbb{R}^{n}\left\vert \Phi ^{\mathbb{Q}}\left(
\mathbf{v,}T\right) \right. \right\} =\Phi ^{\mathbb{Q}}\left( \mathbf{0,}%
T\right) =1.  \label{max111}
\end{equation}%
Let $\nu _{s}\in \left( 0,1\right) ,1\leq s\leq n$. From (\ref{phiqvt}), (%
\ref{asymp 1}) and (\ref{max111}) it follows that $\left\vert \Phi ^{\mathbb{%
Q}}\left( \mathbf{v,}T\right) \right\vert $ can be asymptotically majorated
as $\left\vert v_{s}\right\vert \rightarrow \infty ,1\leq s\leq n,\mathbf{%
v\in }T_{n}^{^{\prime }},$
\[
\left\vert \Phi ^{\mathbb{Q}}\left( \mathbf{v,}T\right) \right\vert \lesssim
\prod\limits_{s=1}^{n}\left\vert \exp \left( -T\psi _{s}^{\mathbb{Q}}\left(
v_{s}\right) \right) \right\vert
\]%
\[
\simeq \prod\limits_{s=1}^{n}\left\vert \exp \left( -T\left( -i\mu
_{s}x_{s}-\Gamma \left( -\nu _{s}\right) \left( c_{s,+}\exp \left( \frac{%
-i\nu _{s}\pi }{2}\right) +c_{s,-}\exp \left( \frac{i\nu _{s}\pi }{2}\right)
\right) \left\vert x_{s}\right\vert ^{\nu }\right) \right) \right\vert
\]%
\[
\lesssim \prod\limits_{s=1}^{n}\exp \left( -T\left( -\Gamma \left( -\nu
_{s}\right) \right) \cos \left( \frac{\nu _{s}\pi }{2}\right) \left(
c_{s,+}+c_{s,-}\right) \left\vert x_{s}\right\vert ^{\nu _{s}}\right)
\]%
\begin{equation}
:=\prod\limits_{s=1}^{n}\exp \left( -Td_{s}\left\vert x_{s}\right\vert ^{\nu
_{s}}\right) =\exp \left( -\sum_{s=1}^{n}\left\vert \left( d_{s}T\right)
^{\nu _{s}^{-1}}x_{s}\right\vert ^{\nu _{s}}\right)  \label{majorant-phi}
\end{equation}%
where
\[
d_{s}:=-\Gamma \left( -\nu _{s}\right) \cos \left( \frac{\nu _{s}\pi }{2}%
\right) \left( c_{s,+}+c_{s,-}\right) ,
\]%
$x_{s}=\Re q
v_{s}$ and $d_{s}>0,1\leq s\leq n.$ For a fixed $0<\varrho
<1 $ consider the set
\[
\Omega \left( \Phi ^{\mathbb{Q}}\left( \mathbf{v,}t\right) ,\varrho \right)
:=\left\{ \mathbf{v}\in \mathbb{R}^{n},\left\vert \Phi ^{\mathbb{Q}}\left(
\mathbf{v,}t\right) \right\vert \geq \varrho \right\} .
\]%
From (\ref{majorant-phi}) it follows that
\[
\Omega \left( \Phi ^{\mathbb{Q}}\left( \mathbf{v,}T\right) ,\varrho \right)
\subset \Omega _{\varrho }^{^{\prime }}:=\left\{ \mathbf{x}\in \mathbb{R}%
^{n},\exp \left( -\sum_{s=1}^{n}\left\vert \left( d_{s}T\right) ^{\nu
_{s}^{-1}}x_{s}\right\vert ^{\nu _{s}}\right) \geq \varrho \right\}
\]%
or%
\[
\Omega _{\varrho }^{^{\prime }}:=\left\{ \mathbf{x}\in \mathbb{R}%
^{n},\sum_{s=1}^{n}\left\vert \left( \frac{d_{s}T}{\ln \varrho ^{-1}}\right)
^{\nu _{s}^{-1}}x_{s}\right\vert ^{\nu _{s}}\leq 1\right\} .
\]%
Now we have
\[
\mathrm{Card}\left( \Omega _{\varrho }^{^{\prime }}\cap \mathbb{Z}%
^{n}\right) \simeq \mathrm{Vol}_{n}\left( \Omega _{\varrho }^{^{\prime
}}\right) ,
\]%
as $\varrho \rightarrow 0$ and
\[
\mathrm{Vol}_{n}\left( \Omega _{\varrho }^{^{\prime }}\right)
=
prod_{s=1}^{n}\left( \frac{\ln \varrho ^{-1}}{d_{s}T}\right) ^{\nu
_{s}^{-1}}B\left( \nu _{1},\cdot \cdot \cdot ,\nu _{n}\right) ,
\]%
where%
\[
B\left( \nu _{1},\cdot \cdot \cdot ,\nu _{n}\right) :=\left\{ \mathbf{x}%
=\left( x_{1},\cdot \cdot \cdot ,x_{n}\right) \in \mathbb{R}%
^{n},\sum_{s=1}^{n}\left\vert x_{s}\right\vert ^{\nu _{s}}\leq 1\right\}
\]%
and $\nu _{1}>0,\cdot \cdot \cdot ,\nu _{n}>0.$ It is known that%
\[
\mathrm{Vol}_{n}B\left( \nu _{1},\cdot \cdot \cdot ,\nu _{n}\right) =2^{n}%
\frac{\prod_{s=1}^{n}\Gamma \left( 1+\nu _{s}\right) }{\Gamma \left(
1+\sum_{s=1}^{n}\nu _{s}\right) }.
\]%
Hence%
\[
\mathrm{Card}\left( \Omega \left( \Phi ^{\mathbb{Q}}\left( \mathbf{v,}%
T\right) ,\varrho \right) \cap \mathbb{Z}^{n}\right)
\]%
\[
\lesssim \varkappa _{n}\left( \ln \varrho ^{-1}\right) ^{\sum_{s=1}^{n}\nu
_{s}^{-1}},\varrho \rightarrow 0.
\]

{\bf Theorem 5.}
\label{riesz theorem} 
\begin{em}
(F. Riesz, \cite{riesz1} p. 102) Let $\omega
_{n}\left( \mathbf{x}\right) $, $n\in \mathbb{N}$ be othonormal system of
functions which is uniformly boned over $\Omega \subset \mathbb{R}^{n}$,
\[
\sup_{\mathbf{x}\in \mathbb{R}^{n}}\left\vert \omega _{n}\left( \mathbf{x}%
\right) \right\vert \leq L.
\]%
If $1\leq p\leq 2$ then there is such $f\in L_{p}$ such that
\[
\left\Vert f\right\Vert _{p^{^{\prime }}}\leq L^{2/p-1}\left\Vert
c\right\Vert _{p},
\]%
where $1/p+1/p^{^{\prime }}=1$ and
\[
\left\Vert c\right\Vert _{p}:=\left( \sum_{k=1}^{\infty }\left\vert
c_{k}\right\vert ^{1/p}\right) ^{1/p},c_{k}=\int_{\Omega }f\overline{\omega }%
_{n}d\mathbf{x.}
\]
\end{em}

{\bf Theorem 6.} \label{0000}
\begin{em}
Let $2\leq p\leq \infty ,1/p+1/p^{^{\prime }}=1,$%
\[
e_{n}:=P^{-n/p^{^{\prime }}}\left( \frac{\varkappa _{n}}{p^{^{\prime }}}%
\sum_{s=1}^{n}\nu _{s}^{-1}\right) ^{1/p^{^{\prime }}},
\]%
\[
\eta _{n}:=e_{n}\varkappa _{n}^{\left( \sum_{s=1}^{n}\nu _{s}^{-1}-1\right)
\left( p^{^{\prime }}\right) ^{-1}}
\]%
and%
\[
M:=\varkappa _{n}\left( \ln R\right) ^{\sum_{s=1}^{n}\nu _{s}^{-1}}.
\]%
Then in our notations%
\[
E\left( M\right) :=\left\Vert \frac{1}{P^{n}}\sum_{\mathbf{m}\in \mathbb{Z}%
^{n}}\Phi ^{\mathbb{Q}}\left( -\frac{\mathbf{m}}{P},T\right) \exp \left(
i\left\langle \frac{2\pi \mathbf{m}}{P},\mathbf{x}\right\rangle \right)
\right.
\]%
\[
\left. -\frac{1}{P^{n}}\sum_{\mathbf{m}\in \mathbb{Z}^{n}\cap \Omega
_{1/R}^{^{\prime }}}\Phi ^{\mathbb{Q}}\left( -\frac{\mathbf{m}}{P},T\right)
\exp \left( i\left\langle \frac{2\pi \mathbf{m}}{P},\mathbf{x}\right\rangle
\right) \right\Vert _{L_{p}\left( \frac{P}{2}Q_{n}\right) }
\]%
\[
\lesssim \eta _{n}\exp \left( -\varkappa _{n}^{-\left( \sum_{s=1}^{n}\nu
_{s}^{-1}\right) ^{-1}}M^{\left( \sum_{s=1}^{n}\nu _{s}^{-1}\right)
^{-1}}\right) M^{\left( 1-\sum_{s=1}^{n}\nu _{s}^{-1}\right) ^{-1}\left(
p^{^{\prime }}\right) ^{-1}},
\]%
as $M\rightarrow \infty $.
\end{em}

{\bf Proof}
Observe that the system of functions%
\[
\chi _{\mathbf{m}}\left( \mathbf{x}\right) :=P^{-n/2}\exp \left(
i\left\langle \frac{2\pi \mathbf{m}}{P},\mathbf{x}\right\rangle \right) ,%
\mathbf{m\in }\mathbb{Z}^{n},\mathbf{x\in }\frac{P}{2}Q_{n}
\]%
is uniformly bounded
\[
\left\vert \chi _{\mathbf{m}}\left( \mathbf{x}\right) \right\vert \leq
P^{-n/2},\forall \mathbf{m\in }\mathbb{Z}^{n}
\]%
and orthonormal in $L_{2}\left( \frac{P}{2}Q_{n}\right) $. Let $\varrho
\rightarrow \infty $ then in our notations%
\[
\mathrm{Vol}_{n}\left( \Omega _{1/\rho }^{^{\prime }}\right) \simeq
\varkappa _{n}\left( \ln \rho \right) ^{\sum_{s=1}^{n}\nu _{s}^{-1}}.
\]
Applying Theorem \ref{riesz theorem} we get
\[
\left\Vert \frac{1}{P^{n}}\sum_{\mathbf{m}\in \left( \mathbb{R}^{n}\setminus
\Omega _{1/R}^{^{\prime }}\right) \cap \mathbb{Z}^{n}}\Phi ^{\mathbb{Q}%
}\left( -\frac{\mathbf{m}}{P},T\right) \exp \left( i\left\langle \frac{2\pi
\mathbf{m}}{P},\mathbf{x}\right\rangle \right) \right\Vert _{L_{p}\left(
\frac{P}{2}Q_{n}\right) }
\]%
\[
=\left\Vert \frac{1}{P^{n/2}}\sum_{\mathbf{m}\in \left( \mathbb{R}%
^{n}\setminus \Omega _{1/R}^{^{\prime }}\right) \cap \mathbb{Z}^{n}}\Phi ^{%
\mathbb{Q}}\left( -\frac{\mathbf{m}}{P},T\right) \chi _{\mathbf{m}}\left(
\mathbf{x}\right) \right\Vert _{L_{p}\left( \frac{P}{2}Q_{n}\right) }
\]%
\[
\lesssim \frac{1}{P^{n/2}}\left( \int_{R}^{\infty }\rho ^{-p^{^{\prime }}}d%
\mathrm{Vol}_{n}\left( \Omega _{1/\rho }^{^{\prime }}\right) \right)
^{1/p^{^{\prime }}}
\]%
\[
\simeq \frac{1}{P^{n/2}}\left( \frac{\varkappa _{n}}{p^{^{\prime }}}%
\sum_{s=1}^{n}\nu _{s}^{-1}\right) ^{1/p^{^{\prime }}}R^{-1}\left( \ln
R\right) ^{\left( \sum_{s=1}^{n}\nu _{s}^{-1}-1\right) \left( p^{^{\prime
}}\right) ^{-1}},R\rightarrow \infty .
\]%
It means that using
\[
M:=\varkappa _{n}\left( \ln R\right) ^{\sum_{s=1}^{n}\nu _{s}^{-1}}
\]%
harmonics from\textrm{\ }$\Omega _{1/R}^{^{\prime }}$ we get the error $%
E\left( M\right) $ of approximation
\[
E\left( M\right) \lesssim \frac{L^{2/p^{^{\prime }}-1}}{P^{n/2}}\left( \frac{%
\varkappa _{n}}{p^{^{\prime }}}\sum_{s=1}^{n}\nu _{s}^{-1}\right)
^{1/p^{^{\prime }}}R^{-1}\left( \ln R\right) ^{\left( \sum_{s=1}^{n}\nu
_{s}^{-1}-1\right) \left( p^{^{\prime }}\right) ^{-1}}
\]%
\[
=P^{-n/p^{^{\prime }}}\left( \frac{\varkappa _{n}}{p^{^{\prime }}}%
\sum_{s=1}^{n}\nu _{s}^{-1}\right) ^{1/p^{^{\prime }}}R^{-1}\left( \ln
R\right) ^{\left( \sum_{s=1}^{n}\nu _{s}^{-1}-1\right) \left( p^{^{\prime
}}\right) ^{-1}}
\]%
\[
:=e_{n}R^{-1}\left( \ln R\right) ^{\left( \sum_{s=1}^{n}\nu
_{s}^{-1}-1\right) \left( p^{^{\prime }}\right) ^{-1}},R\rightarrow \infty ,
\]%
or%
\[
E\left( M\right) \lesssim \eta _{n}\exp \left( -\varkappa _{n}^{-\left(
\sum_{s=1}^{n}\nu _{s}^{-1}\right) ^{-1}}M^{\left( \sum_{s=1}^{n}\nu
_{s}^{-1}\right) ^{-1}}\right) M^{\left( 1-\sum_{s=1}^{n}\nu
_{s}^{-1}\right) ^{-1}\left( p^{^{\prime }}\right) ^{-1}}
\]%
as $M\rightarrow \infty $.

{\bf Corollary 1.}
\label{collorary-conv}
\begin{em}
\[
\left\Vert p_{t}^{\mathbb{Q}}\left( \mathbf{x}\right) -\frac{1}{P^{n}}\sum_{%
\mathbf{m}\in \left( \mathbb{R}^{n}\setminus \Omega _{1/R}^{^{\prime
}}\right) \cap \mathbb{Z}^{n}}\Phi ^{\mathbb{Q}}\left( -\frac{\mathbf{m}}{P}%
,T\right) \exp \left( i\left\langle \frac{2\pi \mathbf{m}}{P},\mathbf{x}%
\right\rangle \right) \right\Vert _{L_{p}\left( \frac{P}{2}Q_{n}\right) }
\]%
\[
\leq \epsilon P^{n/p}+\eta _{n}\exp \left( -\varkappa _{n}^{-\left(
\sum_{s=1}^{n}\nu _{s}^{-1}\right) ^{-1}}M^{\left( \sum_{s=1}^{n}\nu
_{s}^{-1}\right) ^{-1}}\right) M^{\left( 1-\sum_{s=1}^{n}\nu
_{s}^{-1}\right) ^{-1}\left( p^{^{\prime }}\right) ^{-1}},
\]%
as $M\rightarrow \infty $.
\end{em}

\section{The problem of optimal approximation of density functions and $n$%
-widths}

In this section we discuss the problem of optimal approximation of density
functions using a wide range of approximation methods. In problems of
optimal recovery arise quantities which are known as cowidths. Let $\left(
X,d\right) $ be a given metric (Banach) space, $Y$ a certain set (coding
set), $A\subset X$, $\Phi $ a family of mappings $\phi :A\rightarrow Y$,
then the respective cowidth can be defined as
\[
\mathrm{co}^{\Phi }\left( A,X\right) =\inf_{\phi \in \Phi }\sup_{y\in \phi
\left( A\right) }\mathrm{diam}\left\{ \phi ^{-1}\left( y\right) \cap
A\right\} ,
\]%
where%
\[
\phi ^{-1}\left( y\right) =\left\{ x\left\vert x\in X,\phi \left( x\right)
=\phi \left( y\right) \right. \right\} .
\]%
In particular, let $Y$ be $\mathbb{R}^{m}$ and $\Phi :\mathrm{lin}\left(
A\right) \rightarrow \mathbb{R}^{m}$ be a linear application, $\Phi =%
\mathcal{L}\left( \mathrm{lin}\left( A\right) ,\mathbb{R}^{m}\right) $, then
we get a linear cowidth $\lambda ^{m}\left( A,X\right) $. It is easy to
check that $\lambda ^{m}=2d^{m}$, where $d^{m}$ is the Gelfand's $m$-width
defined by
\[
d^{m}\left( A,X\right) =\inf \left\{ L_{-m}\subset X\left\vert
\sup \left\{ x\left\vert x\in A\cap L_{-m}\right. \right\} \right. \right\}
,
\]%
where $\inf $ is taken over all subspaces $L_{-m}$ of codimension $m$.
Letting $Y$ be the set of all $m$-dimensional complexes in $X$ and $\Phi
=C\left( A,Y\right) $ be the set of all continuous mappings $\phi
:A\rightarrow Y$, then we get Alexandrov's cowidths $a^{m}\left( A,X\right) $%
.

Let $P\in \mathbb{N}$ and $\left\{ \varrho _{k}\left( \mathbf{x}\right)
,k\in \mathbb{N}\right\} $ be a set of continuous orthonormal functions on
the $n$-dimensional torus, $P\mathbb{T}^{n}:=\mathbb{R}^{n}/P\mathbb{Z}^{n}.$
Let $\Lambda :=\left\{ \lambda _{k},k\in \mathbb{N}\right\} $ be a fixed
sequence of complex numbers such that $\left\vert \lambda _{1}\right\vert
\geq $ $\left\vert \lambda _{2}\right\vert \geq \cdot \cdot \cdot \geq
\left\vert \lambda _{m}\right\vert \geq \cdot \cdot \cdot $. For any $f\in
L_{1}\left( P\mathbb{T}^{n}\right) $ we can construct a formal Fourier
series
\[
\mathtt{s}\left[ f\right] =\sum_{k=1}^{\infty }c_{k}\left( f\right) \varrho
_{k}\mathbf{,}c_{k}\left( f\right) :=\int_{P\mathbb{T}^{n}}f\varrho _{k}dx.
\]%
Consider set of functions
\[
\Lambda :=\left\{ f\left\vert \left\vert c_{k}\left( f\right) \right\vert
\leq \left\vert \lambda _{k}\right\vert ,k\in \mathbb{N}\right. \right\} .
\]%
It is easy to check that $\Lambda $ is a convex and symmetric set. Also, $%
\Lambda $ is compact in $\ C\left( P\mathbb{T}^{n}\right) $ if
\[
\sum_{k=1}^{\infty }\left\vert \lambda _{k}\right\vert <\infty .
\]
{\bf Theorem 7.}
\label{lower bound} 
\begin{em}
Let $\sum_{k=1}^{\infty }\lambda _{k}<\infty $ then
\[
a^{m}\left( \Lambda ,L_{1}\left( P\mathbb{T}^{n}\right) \right) \geq
2^{-1}\left\vert \lambda _{m+1}\right\vert .
\]
\end{em}

{\bf Proof}
Let us remind some basic definitions. Let $X$ be a Banach space with the
norm $\left\Vert \cdot \right\Vert $ and the unit ball $B$ and $A$ be a
convex, compact, centrally symmetric subset of $X$. Let $L_{m+1}$ be an $%
(m+1)$-dimensional subspace in $X$. Bernstein's $m$-width is defined as
\[
b_{m}\left( A,X\right) =\sup \left\{ L_{m+1}\subset X\left\vert \sup \left\{
\epsilon >0\left\vert \epsilon B\cap L_{m+1}\subset A\right. \right\}
\right. \right\} .
\]%
The Alexandrov $m$-width is the value
\[
a_{m}\left( A,X\right) =\inf_{\Theta _{m}\subset X}\inf_{\sigma
:A\rightarrow \Theta _{m}}\sup \left\{ x\left\vert x\in A,\left\Vert
x-\sigma \left( x\right) \right\Vert \right. \right\} ,
\]%
where the infimum is taken over all $m$-dimensional complexes $\Theta _{m}$,
lying in $X$ and all continuous mappings $\sigma :A\rightarrow \Theta _{m}$.

The Uryson width $u_{m}\left( A,X\right) $ is the infimum of those $\epsilon
>0$ for which there exists a covering of $A$ by open sets (in the sense of
topology induced by the norm $\left\Vert \cdot \right\Vert $ in $X$) of
diameter $<\epsilon $ in $X$ of multiplicity $m+1$ (i.e., such that each
point is covered by $\leq m+1$ sets and some point is covered by exactly $%
m+1 $ sets). Observe that the width $u_{m}\left( A,X\right) $ was introduced
by Uryson \cite{uryson} and inspired by the Lebesgue-Brouwer definition of
dimension.

It is known \cite{tikhomirov1} p.190 that for any compact set $A$ in a
Banach space $X$
\[
b_{m}\left( A,X\right) \leq 2a_{m}\left( A,X\right)
\]%
\cite{tikhomirov2} p.222,%
\[
a_{m}\left( A,X\right) \leq u_{n}\left( A,X\right)
\]%
and \cite{tikhomirov1} p.190,%
\[
u_{m}\left( A,X\right) \leq a^{m}\left( A,X\right) .
\]%
Hence (see \cite{kushpel-tozoni1} for more details)%
\begin{equation}
b_{m}\left( A,X\right) \leq 2a^{m}\left( A,X\right) .  \label{alexandrov1}
\end{equation}%
Let us fix
\[
L_{m+1}=\mathrm{lin}\left\{ \varrho _{k},1\leq k\leq m+1\right\}
\]%
and consider the set
\[
QL_{m+1}:=\left\{ t_{m+1}=\sum_{k=1}^{m+1}c_{k}\varrho _{k},\left\vert
c_{k}\right\vert \leq \left\vert \lambda _{m+1}\right\vert \right\} .
\]%
Clearly%
\[
QL_{m+1}\subset \Lambda
\]%
and%
\[
\left\vert \lambda _{m+1}\right\vert B_{1}\left( P\mathbb{T}^{n}\right) \cap
L_{m+1}\subset QL_{m+1},
\]%
where%
\[
B_{1}\left( P\mathbb{T}^{n}\right) :=\left\{ f\left\vert \left\Vert
f\right\Vert _{1}\leq 1\right. \right\} .
\]%
It means that
\begin{equation}
b_{m}\left( \Lambda ,L_{1}\left( P\mathbb{T}^{n}\right) \right) \geq
\left\vert \lambda _{m+1}\right\vert .  \label{bernst1}
\end{equation}%
Finally, comparing (\ref{bernst1}) and (\ref{alexandrov1}) we get%
\[
a^{m}\left( \Lambda ,L_{1}\left( P\mathbb{T}^{n}\right) \right) \geq
2^{-1}\left\vert \lambda _{m+1}\right\vert .
\]


For simplicity assume that $\mathrm{A}=0$. In this case
\[
\Phi ^{\mathbb{Q}}\left( \mathbf{v,}t\right) =\prod\limits_{s=1}^{n}\exp
\left( -T\psi _{s}^{\mathbb{Q}}\left( v_{s}\right) \right) .
\]%
Consider function class
\[
\Lambda :=\left\{ f\left( \mathbf{x}\right) =\sum_{\mathbf{m}\in \mathbb{Z}%
^{n}}c_{\mathbf{m}}\varrho _{\mathbf{m}}\left( \mathbf{x}\right) \right\} ,
\]%
where $\left\vert c_{\mathbf{m}}\right\vert \leq \lambda _{\mathbf{m}},$%
\[
\varrho _{\mathbf{m}}\left( \mathbf{x}\right) :=\exp \left( i\left\langle
\frac{2\pi \mathbf{m}}{P},\mathbf{x}\right\rangle \right) ,\mathbf{m}\in
\mathbb{Z}^{n}
\]%
and%
\[
\lambda _{\mathbf{m}}=\frac{1}{P^{n/2}}\left\vert \Phi ^{\mathbb{Q}}\left( -%
\frac{\mathbf{m}}{P}\mathbf{,}T\right) \right\vert .
\]%
Using Theorem \ref{lower bound} and Theorem \ref{m-term1} we get the
following statement.

{\bf Corollary 1.}  \label{0002}
\begin{em}
In our notations
\[
a^{M}\left( \Lambda ,L_{1}\left( P\mathbb{T}^{n}\right) \right) \gtrsim
\frac{1}{2}\exp \left( -\left( \varkappa _{n}^{-1}M\right) ^{\left(
\sum_{s=1}^{n}\nu _{s}^{-1}\right) ^{-1}}\right) ,M\rightarrow \infty .
\]
\end{em}

\section{Pricing of basket options}

In applications it is important to construct such pricing theory which
includes a wide range of reward functions $\varphi $. For instance, the
European call reward function which is given by
\[
\varphi =\varphi \left( x_{1},\cdot \cdot \cdot ,x_{n}\right) =\left(
S_{0,1}\exp \left( x_{1}\right) -\sum_{j=2}^{n}S_{0,j}\exp \left(
x_{j}\right) -K\right) _{+}
\]%
admits an exponential grows with respect to $x_{1}$ as $x_{1}\rightarrow
\infty $. Hence we need to introduce the following definition.

{\bf Definition 1.} \label{gavno}
\begin{em}
We say that the reward function $\varphi $ is adopted to
the model process $\mathbf{U}_{t}=\left\{ \mathbf{U}_{t},t\in \mathbb{R}%
_{+}\right\} $\textit{\ if }$\ \mathbb{E}^{\mathbb{Q}}\left[ \varphi \right]
<\infty $\textit.
\end{em}

Clearly, if $\mathbb{E}^{\mathbb{Q}}\left[ \varphi \right] =\infty $ then
the option can not be priced. Remind that the operator of expectation is
taken with respect to the density function $p_{t}^{\mathbb{Q}}$\ which
satisfies the equivalent martingale measure condition (\ref{EMM condition}).

The next statement reduces European call reward function to the canonic form.

{\bf Lemma 1.} 
\label{shift1} 
\begin{em}
In our notations
\[
V=K\exp \left( -rT\right) \int_{\mathbb{R}^{n}}\left( \exp \left(
y_{1}\right) -\sum_{j=2}^{n}\exp \left( y_{j}\right) -1\right) _{+}p_{T}^{%
\mathbb{Q}}\left( \mathbf{y-b}\right) d\mathbf{y,}
\]%
where
\[
\mathbf{b}:=\left( a_{1},\cdot \cdot \cdot ,a_{n}\right) ,b_{j}=\ln \left(
\frac{S_{0,j}}{K}\right) ,1\leq j\leq n.
\]
\end{em}

{\bf Proof}
Remind that $V=\exp \left( -rT\right) \mathbb{E}^{\mathbb{Q}}\left[ \varphi %
\right] $. In our case%
\[
\varphi =\left( S_{1,T}-\sum_{j=2}^{n}S_{j,T}-K\right) _{+},
\]%
where
\[
S_{j,T}=S_{0,j}\exp \left( U_{j,T}\right) ,1\leq j\leq n.
\]%
It means that
\[
V=\exp \left( -rT\right) \int_{\mathbb{R}^{n}}\left( S_{0,1}\exp \left(
x_{1}\right) -\sum_{j=2}^{n}S_{0,j}\exp \left( x_{j}\right) -K\right)
_{+}p_{T}^{\mathbb{Q}}\left( \mathbf{x}\right) d\mathbf{x,}
\]%
\[
=K\exp \left( -rT\right)
\]%
\[
\times \int_{\mathbb{R}^{n}}\left( \exp \left( x_{1}+\ln \left( \frac{S_{0,1}%
}{K}\right) \right) -\sum_{j=2}^{n}\exp \left( x_{j}+\ln \left( \frac{S_{0,j}%
}{K}\right) \right) -1\right) _{+}p_{T}^{\mathbb{Q}}\left( \mathbf{x}\right)
d\mathbf{x,}
\]%
where $S_{0,j},1\leq j\leq n$ are the respective spot prices. Making change
of variables
\[
y_{j}=x_{j}+\ln \left( \frac{S_{0,j}}{K}\right) ,1\leq j\leq n
\]%
we get%
\[
V=K\exp \left( -rT\right) \int_{\mathbb{R}^{n}}\left( \exp \left(
y_{1}\right) -\sum_{j=2}^{n}\exp \left( y_{j}\right) -1\right) _{+}p_{T}^{%
\mathbb{Q}}\left( \mathbf{y-b}\right) d\mathbf{y,}
\]%
where
\[
\mathbf{b}:=\left( b_{1},\cdot \cdot \cdot ,b_{n}\right) ,b_{j}=\ln \left(
\frac{S_{0,j}}{K}\right) ,1\leq j\leq n.
\]

We will need the following result \cite{hurd}.

{\bf Theorem 8.}
\label{hurd-result} 
\begin{em}
Let $n\geq 2$. For any real numbers $\mathbf{\epsilon }%
=\left( \epsilon _{1},\cdot \cdot \cdot ,\epsilon _{n}\right) $ with $%
\epsilon _{m}>0$ for $2\leq m\leq n$ and $\epsilon
_{1}<-1-\sum_{m=2}^{n}\epsilon _{m},$%
\[
\left( \exp \left( x_{1}\right) -\sum_{m=2}^{n}\exp \left( x_{m}\right)
-1\right) _{+}=\left( 2\pi \right) ^{-n}\int_{\mathbb{R}^{n}+i\mathbf{%
\epsilon }}\exp \left( i\left\langle \mathbf{u},\mathbf{x}\right\rangle
\right) \mathbf{F}S\left( \mathbf{u}\right) d\mathbf{u,}
\]%
where $\mathbf{x}=\left( x_{1},\cdot \cdot \cdot ,x_{n}\right) $ and for $%
\mathbf{u=}\left( u_{1},\cdot \cdot \cdot ,u_{n}\right) \in \mathbb{C}^{n}$%
\[
\mathbf{F}S\left( \mathbf{u}\right) =\frac{\Gamma \left( i\left(
u_{1}+\sum_{m=2}^{n}u_{m}\right) -1\right) \prod\limits_{m=2}^{n}\Gamma
\left( -iu_{m}\right) }{\Gamma \left( iu_{1}+1\right) }.
\]
\end{em}

The next statement gives a general approximation formula for the European
call options which is important in various applications. Observe that it
does not show the rate of convergence. This problem will be discussed later,
it explains just how to construct the approximation formula.

{\bf Theorem 9.} \label{approximant 1} 
\begin{em}
Let%
\[
\mathbf{b}:=\left( b_{1},\cdot \cdot \cdot ,b_{n}\right) ,b_{j}=\ln \left(
\frac{S_{0,j}}{K}\right) ,1\leq j\leq n
\]%
and $\mathbf{\epsilon }=\left( \epsilon _{1},\cdot \cdot \cdot ,\epsilon
_{n}\right) ,-i\mathbf{\epsilon }\in T^{^{\prime }},$ $2\leq m\leq
n,\epsilon _{1}<-1-\sum_{m=2}^{n}\epsilon _{m}.$ Then formal approximation
formula can be written as%
\[
V\approx \frac{K\exp \left( -rT-\left\langle \mathbf{b,\epsilon }%
\right\rangle \right) }{P^{n}}\sum_{\mathbf{m}\in \Omega _{1/R}^{^{\prime
}}}\Phi ^{\mathbb{Q}}\left( \frac{\mathbf{m}}{P}+i\mathbf{\epsilon }%
,T\right) \exp \left( i\left\langle \frac{2\pi \mathbf{m}}{P},\mathbf{b}%
\right\rangle \right)
\]%
\[
\times \frac{\Gamma \left( -\frac{2\pi i}{P}\left(
m_{1}+\sum_{s=2}^{n}m_{s}\right) +1\right) \prod\limits_{s=2}^{n}\Gamma
\left( \frac{2\pi i}{P}m_{s}\right) }{\Gamma \left( 1-\frac{2\pi im_{1}}{P}%
\right) },
\]
where $R\rightarrow \infty ,P\rightarrow \infty $ and
\[
\Omega _{1/R}^{^{\prime }}=\left\{ \mathbf{x}\in \mathbb{R}%
^{n},\sum_{s=1}^{n}\left\vert \left( \frac{d_{s}T}{\ln R}\right) ^{\nu
_{s}^{-1}}x_{s}\right\vert ^{\nu _{s}}\leq 1\right\}
\]%
and%
\[
d_{s}=-\Gamma \left( -\nu _{s}\right) \cos \left( \frac{\nu _{s}\pi }{2}%
\right) \left( c_{s,+}+c_{s,-}\right) ,0<\nu _{s}<1.
\]
\end{em}

{\bf Proof}
Applying Lemma \ref{shift1} we get%
\[
V=\exp \left( -rT\right) \mathbb{E}^{\mathbb{Q}}\left[ \varphi \right]
\]%
\[
=\exp \left( -rT\right) \int_{\mathbb{R}^{n}}\left( S_{0,1}\exp \left(
x_{1}\right) -\sum_{j=2}^{n}S_{0,j}\exp \left( x_{j}\right) -K\right)
_{+}p_{T}^{\mathbb{Q}}\left( \mathbf{x}\right) d\mathbf{x,}
\]%
\[
=K\exp \left( -rT\right) \int_{\mathbb{R}^{n}}\left( \exp \left(
y_{1}\right) -\sum_{j=2}^{n}\exp \left( y_{j}\right) -1\right) _{+}p_{T}^{%
\mathbb{Q}}\left( \mathbf{y-b}\right) d\mathbf{y,}
\]%
where%
\[
\mathbf{b}:=\left( b_{1},\cdot \cdot \cdot ,b_{n}\right) ,b_{j}=\ln \left(
\frac{S_{0,j}}{K}\right) ,1\leq j\leq n.
\]%
Assume that $-i\mathbf{\epsilon }\in T_{n}^{^{\prime }}$. Applying Cauchy's
theorem $n$ times in the tube $T_{n}^{^{\prime }}$, which is justified by (%
\ref{asymp 1}), we get%
\[
p_{T}^{\mathbb{Q}}\left( \mathbf{y}\right) =\left( 2\pi \right) ^{-n}\int_{%
\mathbb{R}^{n}}\exp \left( -i\left\langle \mathbf{y},\mathbf{x}\right\rangle
\right) \Phi ^{\mathbb{Q}}\left( \mathbf{x},T\right) d\mathbf{x}
\]%
\[
=\left( 2\pi \right) ^{-n}\int_{\mathbb{R}^{n}-i\mathbf{\epsilon }}\exp
\left( -i\left\langle \mathbf{y},\mathbf{x}\right\rangle \right) \Phi ^{%
\mathbb{Q}}\left( \mathbf{x},T\right) d\mathbf{x}
\]%
\[
=\left( 2\pi \right) ^{-n}\int_{\mathbb{R}^{n}}\exp \left( -i\left\langle
\mathbf{y},\mathbf{x}+i\mathbf{\epsilon }\right\rangle \right) \Phi ^{%
\mathbb{Q}}\left( \mathbf{x}+i\mathbf{\epsilon },T\right) d\mathbf{x}
\]%
\[
=\exp \left( \left\langle \mathbf{y},\mathbf{\epsilon }\right\rangle \right)
\left( 2\pi \right) ^{-n}\int_{\mathbb{R}^{n}}\exp \left( -i\left\langle
\mathbf{y},\mathbf{x}\right\rangle \right) \Phi ^{\mathbb{Q}}\left( \mathbf{x%
}+i\mathbf{\epsilon },T\right) d\mathbf{x}
\]%
Let $\mathbf{y}\in \frac{P}{2}Q_{n}$ then from Corollary \ref{collorary-conv}
we get
\[
p_{T}^{\mathbb{Q}}\left( \mathbf{y}\right) \approx \exp \left( \left\langle
\mathbf{y},\mathbf{\epsilon }\right\rangle \right) \left( \frac{1}{P^{n}}%
\sum_{\mathbf{m}\in \mathbb{Z}^{n}}\Phi ^{\mathbb{Q}}\left( \mathbf{-}\frac{%
\mathbf{m}}{P}+i\mathbf{\epsilon },T\right) \exp \left( i\left\langle \frac{%
2\pi \mathbf{m}}{P},\mathbf{y}\right\rangle \right) \right)
\]%
and
\[
p_{T}^{\mathbb{Q}}\left( \mathbf{y-b}\right) \approx \exp \left(
\left\langle \mathbf{y-b},\mathbf{\epsilon }\right\rangle \right)
\]%
\[
\times \frac{1}{P^{n}}\sum_{\mathbf{m}\in \Omega _{1/R}^{^{\prime }}}\left(
\Phi ^{\mathbb{Q}}\left( \mathbf{-}\frac{\mathbf{m}}{P}+i\mathbf{\epsilon }%
,T\right) \exp \left( i\left\langle -\frac{2\pi \mathbf{m}}{P},\mathbf{b}%
\right\rangle \right) \right) \exp \left( i\left\langle \frac{2\pi \mathbf{m}%
}{P},\mathbf{y}\right\rangle \right)
\]
Finally, assuming $\mathbf{\epsilon }=\left( \epsilon _{1},\cdot \cdot \cdot
,\epsilon _{n}\right) ,$ $2\leq j\leq n,\epsilon
_{1}<-1-\sum_{j=2}^{n}\epsilon _{j}$, using Theorem \ref{hurd-result} and
the fact that the domain $\Omega _{1/R}^{^{\prime }}$ is centrally symmetric
we obtain
\[
V=K\exp \left( -rT\right) \int_{\mathbb{R}^{n}}\left( \exp \left(
y_{1}\right) -\sum_{j=2}^{n}\exp \left( y_{j}\right) -1\right) _{+}p_{T}^{%
\mathbb{Q}}\left( \mathbf{y-b}\right) d\mathbf{y}
\]%
\[
\approx \frac{K\exp \left( -rT\right) }{P^{n}}\sum_{\mathbf{m}\in \Omega
_{1/R}^{^{\prime }}}\left( \Phi ^{\mathbb{Q}}\left( \mathbf{-}\frac{\mathbf{m%
}}{P}+i\mathbf{\epsilon },T\right) \exp \left( i\left\langle -\frac{2\pi
\mathbf{m}}{P},\mathbf{b}\right\rangle \right) \right)
\]%
\[
\times \int_{\mathbb{R}^{n}}\left( \exp \left( y_{1}\right)
-\sum_{j=2}^{n}\exp \left( y_{j}\right) -1\right) _{+}\exp \left(
\left\langle \mathbf{y-b},\mathbf{\epsilon }\right\rangle \right) \exp
\left( i\left\langle \frac{2\pi \mathbf{m}}{P},\mathbf{y}\right\rangle
\right) d\mathbf{y}
\]%
\[
=\frac{K\exp \left( -rT-\left\langle \mathbf{b,\epsilon }\right\rangle
\right) }{P^{n}}\sum_{\mathbf{m}\in \Omega _{1/R}^{^{\prime }}}\left( \Phi ^{%
\mathbb{Q}}\left( \mathbf{-}\frac{\mathbf{m}}{P}+i\mathbf{\epsilon }%
,T\right) \exp \left( i\left\langle -\frac{2\pi \mathbf{m}}{P},\mathbf{b}%
\right\rangle \right) \right)
\]%
\[
\times \int_{\mathbb{R}^{n}}\left( \left( \exp \left( y_{1}\right)
-\sum_{j=2}^{n}\exp \left( y_{j}\right) -1\right) _{+}\exp \left(
\left\langle \mathbf{y},\mathbf{\epsilon }\right\rangle \right) \right) \exp
\left( i\left\langle \frac{2\pi \mathbf{m}}{P},\mathbf{y}\right\rangle
\right) d\mathbf{y}
\]%
\[
=\frac{K\exp \left( -rT-\left\langle \mathbf{b,\epsilon }\right\rangle
\right) }{P^{n}}\sum_{\mathbf{m}\in \Omega _{1/R}^{^{\prime }}}\left( \Phi ^{%
\mathbb{Q}}\left( \frac{\mathbf{m}}{P}+i\mathbf{\epsilon },T\right) \exp
\left( i\left\langle \frac{2\pi \mathbf{m}}{P},\mathbf{b}\right\rangle
\right) \right)
\]%
\[
\times \frac{\Gamma \left( -\frac{2\pi }{P}\left(
m_{1}+\sum_{s=2}^{n}m_{s}\right) +1\right) \prod\limits_{s=2}^{n}\Gamma
\left( \frac{2\pi i}{P}m_{s}\right) }{\Gamma \left( 1-\frac{2\pi im_{1}}{P}%
\right) }.
\]

The following statement is trivial.

{\bf Lemma 2.} 
\label{trivial1}
\begin{em}
\[
\left( \exp \left( y_{1}\right) -\sum_{j=2}^{n}\exp \left( y_{j}\right)
-1\right) _{+}\leq \min \left\{ \exp \left( y_{1}\right) ,\exp \left(
y_{1}-\sum_{j=2}^{n}\exp \left( y_{j}\right) \right) \right\} .
\]
\end{em}

{\bf Theorem 10.} \label{jopa}
\begin{em}
Let in our notations%
\[
M\left( P,R\right)
\]%
\[
:=\left\Vert \frac{1}{P^{n}}\sum_{\mathbf{m}\in \Omega _{1/R}^{^{\prime
}}}\left( \Phi ^{\mathbb{Q}}\left( \mathbf{-}\frac{\mathbf{m}}{P}+i\mathbf{%
\epsilon },T\right) \exp \left( i\left\langle -\frac{2\pi \mathbf{m}}{P},%
\mathbf{b}\right\rangle \right) \right) \exp \left( i\left\langle \frac{2\pi
\mathbf{m}}{P},\mathbf{y}\right\rangle \right) \right\Vert _{L_{\infty
}\left( \mathbb{R}^{n}\right) },
\]%
\[
L_{\mathbf{\epsilon }}:=\left\Vert \left( \exp \left( y_{1}\right)
-\sum_{j=2}^{n}\exp \left( y_{j}\right) -1\right) _{+}\exp \left(
\left\langle \mathbf{y},\mathbf{\epsilon }\right\rangle \right) \right\Vert
_{L_{1}\left( \mathbb{R}^{n}\right) }
\]%
and $\widetilde{V}$ be the approximant for $V$ from the Theorem \ref%
{approximant 1} then
\[
\delta :=\left\vert V-\widetilde{V}\right\vert
\]%
\[
\leq L_{\mathbf{\epsilon }}\left( \epsilon +\eta _{n}\exp \left( -\varkappa
_{n}^{-\left( \sum_{s=1}^{n}\nu _{s}^{-1}\right) ^{-1}}M^{\left(
\sum_{s=1}^{n}\nu _{s}^{-1}\right) ^{-1}}\right) M^{\left(
1-\sum_{s=1}^{n}\nu _{s}^{-1}\right) ^{-1}}\right)
\]%
\[
+M\left( P,R\right) \int_{\mathbb{R}^{n}\setminus \left( \frac{P}{2}%
-\left\Vert \mathbf{b}\right\Vert _{\infty }\right) Q_{n}}\min \left\{ \exp
\left( y_{1}\right) ,\exp \left( y_{1}-\sum_{j=2}^{n}\exp \left(
y_{j}\right) \right) \right\} \exp \left( \left\langle \mathbf{y},\mathbf{%
\epsilon }\right\rangle \right) .
\]
\end{em}

{\bf Proof} 
Let $\widetilde{V}$ be the approximant for $V$, then assuming that $-i%
\mathbf{\epsilon }\in T^{^{\prime }}$ we get
\[
\delta =\left\vert V-\widetilde{V}\right\vert
\]%
\[
=\left\vert \int_{\mathbb{R}^{n}}\left( \left( \exp \left( y_{1}\right)
-\sum_{j=2}^{n}\exp \left( y_{j}\right) -1\right) _{+}\exp \left(
\left\langle \mathbf{y},\mathbf{\epsilon }\right\rangle \right) \right)
\right.
\]%
\[
\times \left( \left( \int_{\mathbb{R}^{n}}\exp \left( i\left\langle \mathbf{%
y-b,x}\right\rangle \right) \Phi ^{\mathbb{Q}}\left( \mathbf{-x}+i\mathbf{%
\epsilon },T\right) d\mathbf{x}\right) \right.
\]%
\[
\left. \left. -\frac{1}{P^{n}}\sum_{\mathbf{m}\in \Omega _{1/R}^{^{\prime
}}}\left( \Phi ^{\mathbb{Q}}\left( \mathbf{-}\frac{\mathbf{m}}{P}+i\mathbf{%
\epsilon },T\right) \exp \left( i\left\langle -\frac{2\pi \mathbf{m}}{P},%
\mathbf{b}\right\rangle \right) \right) \exp \left( i\left\langle \frac{2\pi
\mathbf{m}}{P},\mathbf{y}\right\rangle \right) \right) d\mathbf{y}%
\right\vert
\]%
\[
:=\int_{\mathbb{R}^{n}}\left( \exp \left( y_{1}\right) -\sum_{j=2}^{n}\exp
\left( y_{j}\right) -1\right) _{+}\exp \left( \left\langle \mathbf{y},%
\mathbf{\epsilon }\right\rangle \right) \mu \left( \mathbf{y}\right) d%
\mathbf{y.}
\]%
From the Theorem \ref{desity 111} and Theorem \ref{approximant 1} it follows
that for a chosen $P>0$ and $M>0$ we have
\[
\delta _{1}:=\left\vert \int_{\mathbb{R}^{n}}\exp \left( i\left\langle
\mathbf{y-b,x}\right\rangle \right) \Phi ^{\mathbb{Q}}\left( \mathbf{-x}+i%
\mathbf{\epsilon },T\right) d\mathbf{x}\right.
\]%
\[
\left. \mathbf{-}\frac{1}{P^{n}}\sum_{\mathbf{m}\in \Omega _{1/R}^{^{\prime
}}}\left( \Phi ^{\mathbb{Q}}\left( \mathbf{-}\frac{\mathbf{m}}{P}+i\mathbf{%
\epsilon },T\right) \exp \left( i\left\langle -\frac{2\pi \mathbf{m}}{P},%
\mathbf{b}\right\rangle \right) \right) \exp \left( i\left\langle \frac{2\pi
\mathbf{m}}{P},\mathbf{y}\right\rangle \right) \right\vert
\]%
\[
\leq \epsilon +\eta _{n}\exp \left( -\varkappa _{n}^{-\left(
\sum_{s=1}^{n}\nu _{s}^{-1}\right) ^{-1}}M^{\left( \sum_{s=1}^{n}\nu
_{s}^{-1}\right) ^{-1}}\right) M^{\left( 1-\sum_{s=1}^{n}\nu
_{s}^{-1}\right) ^{-1}}
\]%
for $\forall \mathbf{y}\in \frac{P}{2}Q_{n}-\mathbf{b}$. Observe that $%
p^{^{\prime }}=1$ in our case. Clearly,
\[
\left( \exp \left( y_{1}\right) -\sum_{j=2}^{n}\exp \left( y_{j}\right)
-1\right) _{+}\exp \left( \left\langle \mathbf{y},\mathbf{\epsilon }%
\right\rangle \right) \in L_{1}\left( \mathbb{R}^{n}\right)
\]%
for a chosen $\mathbf{\epsilon }$. Let%
\[
L_{\mathbf{\epsilon }}=\left\Vert \left( \exp \left( y_{1}\right)
-\sum_{j=2}^{n}\exp \left( y_{j}\right) -1\right) _{+}\exp \left(
\left\langle \mathbf{y},\mathbf{\epsilon }\right\rangle \right) \right\Vert
_{L_{1}\left( \mathbb{R}^{n}\right) }.
\]%
Hence $\frac{P}{2}Q_{n}-\mathbf{b\subset }\left( \frac{P}{2}-\left\Vert
\mathbf{b}\right\Vert _{\infty }\right) Q_{n}$, where $\left\Vert \mathbf{b}%
\right\Vert _{\infty }:=\max \left\{ \left\vert b_{k}\right\vert ,1\leq
k\leq n\right\} $ and
\[
\int_{\left( \frac{P}{2}-\left\Vert \mathbf{b}\right\Vert _{\infty }\right)
Q_{n}}\left( \exp \left( y_{1}\right) -\sum_{j=2}^{n}\exp \left(
y_{j}\right) -1\right) _{+}\exp \left( \left\langle \mathbf{y},\mathbf{%
\epsilon }\right\rangle \right) \mu \left( \mathbf{y}\right) d\mathbf{y}
\]%
\[
\leq L_{\mathbf{\epsilon }}\left( \epsilon +\eta _{n}\exp \left( -\varkappa
_{n}^{-\left( \sum_{s=1}^{n}\nu _{s}^{-1}\right) ^{-1}}M^{\left(
\sum_{s=1}^{n}\nu _{s}^{-1}\right) ^{-1}}\right) M^{\left(
1-\sum_{s=1}^{n}\nu _{s}^{-1}\right) ^{-1}}\right) .
\]%
Then we have
\[
\epsilon ^{^{\prime }}:=\int_{\mathbb{R}^{n}\setminus \left( \frac{P}{2}%
-\left\Vert \mathbf{b}\right\Vert _{\infty }\right) Q_{n}}\left( \exp \left(
y_{1}\right) -\sum_{j=2}^{n}\exp \left( y_{j}\right) -1\right) _{+}\exp
\left( \left\langle \mathbf{y},\mathbf{\epsilon }\right\rangle \right) \mu
\left( \mathbf{y}\right) d\mathbf{y}
\]%
\[
\leq M\left( P,R\right) \left\Vert \left( \exp \left( y_{1}\right)
-\sum_{j=2}^{n}\exp \left( y_{j}\right) -1\right) _{+}\exp \left(
\left\langle \mathbf{y},\mathbf{\epsilon }\right\rangle \right) \right\Vert
_{L_{1}\left( \mathbb{R}^{n}\setminus \left( \frac{P}{2}-\left\Vert \mathbf{b%
}\right\Vert _{\infty }\right) Q_{n}\right) }
\]%
and by Lemma \ref{trivial1}
\[
\epsilon ^{^{\prime }}\leq M\left( P,R\right) \int_{\mathbb{R}^{n}\setminus
\left( \frac{P}{2}-\left\Vert \mathbf{b}\right\Vert _{\infty }\right)
Q_{n}}\left( \exp \left( y_{1}\right) -\sum_{j=2}^{n}\exp \left(
y_{j}\right) -1\right) _{+}\exp \left( \left\langle \mathbf{y},\mathbf{%
\epsilon }\right\rangle \right)
\]%
\[
\leq M\left( P,R\right) \int_{\mathbb{R}^{n}\setminus \left( \frac{P}{2}%
-\left\Vert \mathbf{b}\right\Vert _{\infty }\right) Q_{n}}\min \left\{ \exp
\left( y_{1}\right) ,\exp \left( y_{1}-\sum_{j=2}^{n}\exp \left(
y_{j}\right) \right) \right\} \exp \left( \left\langle \mathbf{y},\mathbf{%
\epsilon }\right\rangle \right) .
\]


\end{document}